# On a Deterministic Property of the Category of *k*-almost Primes: A Deterministic Structure Based on a Linear Function for Redefining the *k*-almost Primes (∃*n*∈ℕ, 1 ≤ *k* ≤ *n*) in Certain Intervals

**Ramin Zahedi**

Logic and Philosophy of Science Research Group, Hokkaido University, Japan.
E-mail: zahedi@let.hokudai.ac.jp
(zahedi.r@gmail.com)

In this paper based on a sort of linear function, a deterministic and simple algorithm with an algebraic structure is presented for calculating all (and only) *k*-almost primes (where ∃*n*∈ℕ, 1 ≤ *k* ≤ *n*) in certain intervals. A theorem has been proven showing a new deterministic property of the category of the *k*-almost primes. Through a linear function that we obtain, an equivalent redefinition of the *k*-almost primes with an algebraic characteristic is identified. Moreover, as an outcome of our function's property some relations which contain new information about the *k*-almost primes (including primes) are presented.

## 1. Introduction

**L**et at first review a short and general summary regarding importance and huge scientific applications of the *k*-almost primes (in particular primes as their base). Indeed, today there are many applications for primes in many scientific fields such as physics, computer science, engineering, security science and chemistry, etc. [6-24, 27, 31-33].

**T**he Quantum mechanical potential is one of the most central concepts in modern quantum deterministic theories. Several studies have shown the relation of algebraic theories, and also primes to the quantum potential [21, 41]. Nowadays, there is interesting speculation that the zeros of the Riemann Zeta function are connected to the energy levels of complex quantum systems [34]. Physicists have reinterpreted the Riemann Zeta function as a (thermodynamic) partition function by defining an abstract numerical 'gas' using the prime numbers [10]; in addition as a result of studying nonlinear dynamics and chaos, they also have discovered at least two instances of fractality within the distribution of prime numbers [6-8].

**C**omplexity theory is a field in theoretical computer science, which attempts to quantify the difficulty of computational tasks and tends to aim at generality while doing so. "Complexity" is measured by various natural computing resources, such as the amount of memory needed, communication bandwidth, time of execution, etc. By analyzing several candidate algorithms for a problem, a most efficient one can be easily identified; for the problem of determining the primality of an integer, the resource that could be examined, is the time of execution [13]. In addition, there are dozens of algorithms in computer science that depend heavily on prime numbers- hashing schemes, sorting schemes, and so on.

**B**ecause of these comprehensive facts regarding importance of primes, it seems that any new research and work about them might find remarkable applications in different scientific fields.

The category of *k*-almost prime numbers is the classified category of both prime and composite numbers. A *k*-almost prime number is a number that has exactly *k* prime factors, where the factors need not be distinct, hence powers of primes are also included [1]. Some recent literatures have proposed calling these numbers primes, biprimes, triprimes, and so on [35]. In this paper we prove the following theorem as a new algebraic property of the category of *k*-almost prime numbers (the initial preprints relating to this paper have been presented in [36]):

**Theorem 1.** Suppose $p_1, p_2, p_3, ..., p_r, p_{r+1}$ are given prime numbers where $p_i$ is the $i^{th}$ prime number, and let $x'_i$ and $t'_i$ are some integer solution of linear equation:

$$p_i x_i - t_i \prod_{k=1}^{i-1} p_k = 1 \qquad (A-1)$$

where $i = 2, 3, 4, ..., r$, then the linear function $Z_r(t, h_i)$:

$$Z_r(t, h_i) = t \prod_{l=1}^{r} p_l + (p_r x'_r - 1)(h_r - 1) + \sum_{j=2}^{r-1}[(p_j x'_j - 1)(h_j - 1) \prod_{q=j+1}^{r} p_q x'_q] - 1 \qquad (A-2)$$

redefines all *k*-almost prime numbers (and only *k*-almost prime numbers) in interval:

$$p_{r+1}^n \leq Z_r(t, h_i) < p_{r+1}^{n+1} \qquad (A-3)$$

where *t* is an integer variable: $t \in \mathbb{Z}$, $h_i \in \{1,2,3,4,...,p_i -1\}$, $\exists n \in \mathbb{N}: 1 \leq k \leq n$, and $r = 3,4,5, ...$; for $r = 1, 2$ we may have

$$Z_1(t, h_i) = 2t + 1, \quad Z_2(t, h_i) = 6t + 4h_2 + 3 \qquad (A-4)$$

(A-4) is obtained from the same way that formula (A-2) is obtained. We will show that in fact function $Z_r(t, h_i)$ defines all integer numbers which are not divisible by $p_1, p_2, ..., p_r$. $x'_i$ in formula (A-2) doesn't have unique values (because equation (A-1) has infinite integer solutions), hence function $Z_r(t, h_i)$ can be written with different coefficients. Using above algebraic property of function $Z_r(t, h_i)$, it is clear that we may also construct a deterministic algorithm for calculating the *k*-almost prime numbers in interval (A-3) and the next proceeding intervals. In fact because the greatest common divisor of $p_i$ and $\prod_{k=1}^{i-1} p_k$ is 1, the linear equation (1) has infinite integer solutions. Equation (1) could be solved not only in a polynomial time algorithm, but also in strongly polynomial time algorithm [37, 38]. In addition, (A-3) also as a linear inequality could be solved in a polynomial time algorithm [39].

**In** a particular case when $n = 1$, function $Z_r(t, h_i)$ **redefines primes** in interval:

$$p_{r+1} \leq Z_r(t, h_i) < p_{r+1}^2 \qquad (A\text{-}5)$$

**In** the next section we prove Theorem 1, and then as some outcome of $Z_r(t, h_i)$ property, some equalities which contain new information about the *k*-almost prime numbers (including primes), are presented.

## 2. Proof of Theorem 1

**It** follows from the definition of primes that in interval $[p_{r+1}^n, p_{r+1}^{n+1})$ any number that is not divisible by any of $p_1, p_2, ..., p_r$ is a *k*-almost prime number (where $\exists n \in \mathbb{N}, 1 \leq k \leq n$) and vice versa. Let $H_r$ be a set of natural numbers $\mathbb{N}$ that are not divisible by $p_1, p_2, ..., p_r$; that is,

$$H_r = \{Z_r \mid Z_r \in \mathbb{N} \text{ and } Z_r \text{ is not divisible by any of } p_1, p_2, ..., p_r\} \qquad (1)$$

Let $E_i$ be the set of natural numbers $\mathbb{N}$ excluding the set of all multipliers of the $i^{th}$ prime number $p_i$; we define sets $E_{i1}, E_{i2}, ..., E_{i(p_i-1)}$:

$$E_{i1} = \{m_{i1} \mid m_{i1} = p_i x_i - 1, \; x_i \in \mathbb{N}\},$$
$$E_{i2} = \{m_{i2} \mid m_{i2} = p_i x_i - 2, \; x_i \in \mathbb{N}\},$$
$$E_{i3} = \{m_{i3} \mid m_{i3} = p_i x_i - 3, \; x_i \in \mathbb{N}\},$$
$$\vdots$$
$$E_{i(p_i-1)} = \{m_{i(p_i-1)} \mid m_{i(p_i-1)} = p_i x_i - (p_i - 1), \; x_i \in \mathbb{N}\} \qquad (2)$$

then

$$E_i = \bigcup_{j=1}^{p_i-1} E_{ij} \qquad (3)$$

which it is equivalent to

$$E_i = \{m_i \mid m_i = p_i x_i - h_i, \; x_i \in \mathbb{N}\} \qquad (4)$$

where $h_i \in \{1, 2, 3, 4, ..., p_i - 1\}$. From the definitions above it follows that for any set $E_{ik}$

and and any set $E_{il}$ ($k \neq l$):

$$E_{ik} \cap E_{il} = \varnothing \tag{5}$$

where $i = 1,2,...,r$ and $k,l = 1,2,...p_i - 1$. It follows from (1), (3) and (4) that

$$H_r = \bigcap_{i=1}^{r} \bigcup_{j=1}^{p_i-1} E_{ij} = \bigcap_{i=1}^{r} E_i \tag{6}$$

Now the following system of linear equations, obtained from (1), (4) and (6) and taking into account relation (5), can define function $Z_r$ in natural numbers:

$$Z_r = p_1 x_1 - h_1 = p_2 x_2 - h_2 = p_3 x_3 - h_3 = .... = p_r x_r - h_r \tag{7}$$

The linear equations in (7) can be re-written as

$$\begin{cases} p_1 x_1 - h_1 = p_2 x_2 - h_2 \\ p_1 x_1 - h_1 = p_3 x_3 - h_3 \\ p_1 x_1 - h_1 = p_4 x_4 - h_4 \\ \quad . \\ \quad . \\ \quad . \\ p_1 x_1 - h_1 = p_r x_r - h_r \end{cases} \tag{8}$$

where $h_i \in \{1,2,3,4,...,p_i - 1\}$ and $i = 1,2,...,r$. Before solving (8), let consider a basic general linear equation in integer numbers $\mathbb{Z}$

$$ax - by = c \tag{9}$$

where x and y are unknown integer values and *a, b,* and *c* are some known integers, and $\text{GCD}(a,b) = 1$, $a > 0, b > 0$. With these conditions, equation (9) has infinite number of integer solutions (positive and negative) and in general

$$x = c\hat{x}' + bt, \quad y = c\hat{y}' + at \tag{10}$$

where $\hat{x}'$, $\hat{y}'$ are some given solution of equation $ax' - by' = 1$ (these given solutions always exist) and $t$ is a free integer parameter $t \in \mathbb{Z}$ [2, 30]. We solve system of equations (8) step by step. Using formula (10) for the first equation in (8) we get

$$p_2 x_2 - p_1 x_1 = h_2 - h_1 = h_2 - 1,$$

with the general solution

$$x_2 = (h_2 - 1)\hat{x}'_2 + p_1 t_1, \quad x_1 = (h_2 - 1)\hat{x}'_1 + p_2 t_1 \quad (11)$$

where $\hat{x}'_1$ and $\hat{x}'_2$ are some given solution of equation $p_2 x'_2 - p_1 x'_1 = 1$ and $t_1$ is a free integer parameter. Using formulas (11) and the second equation of (8) we get

$$p_3 x_3 - p_1 p_2 t_1 = (h_2 - 1) p_1 \hat{x}'_1 + h_3 - 1 \quad (12)$$

with the general solution of

$$x_3 = [h_3 + (h_2 - 1) p_1 \hat{x}'_1 - 1]\hat{x}'_3 + p_1 p_2 t_2,$$
$$t_1 = [h_3 + (h_2 - 1) p_1 \hat{x}'_1 - 1]\hat{t}'_1 + p_3 t_2 \quad (13)$$

where $\hat{x}'_3$ and $\hat{t}'_1$ are some given solution of $p_3 x'_3 - p_1 p_2 t'_1 = 1$ and $t_2$ is any integer value. Using (11), (13) and the third equation of (8) we obtain

$$p_4 x_4 - p_1 p_2 p_3 t_2 = [(h_4 - 1) + (h_3 - 1) p_1 p_2 \hat{t}'_1 + (h_2 - 1) p_1 \hat{x}'_1 p_3 \hat{x}'_3] \quad (14)$$

The general solution of (14) is

$$x_4 = [(h_4 - 1) + (h_3 - 1) p_1 p_2 \hat{t}'_1 + (h_2 - 1) p_1 \hat{x}'_1 p_3 \hat{x}'_3]\hat{x}'_4 + p_1 p_2 p_3 t_3,$$
$$t_2 = [(h_4 - 1) + (h_3 - 1) p_1 p_2 \hat{t}'_1 + (h_2 - 1) p_1 \hat{x}'_1 p_3 \hat{x}'_3]\hat{t}'_2 + p_4 t_3 \quad (15)$$

Continuing this procedure for the proceeding equations, the following general forms of the solutions are obtained:

$$x_i = \hat{x}'_i [(h_i - 1) + (h_{i-1} - 1)\hat{t}'_{i-3} \prod_{l=1}^{i-2} p_l + \sum_{j=2}^{i-2} [(h_j - 1)\hat{t}'_{j-2} \prod_{l=1}^{j-1} p_l \prod_{q=j+1}^{i-1} p_q \hat{x}'_q]] + t_{i-1} \prod_{j=1}^{i-1} p_j,$$

$$t_{i-2} = \hat{t}'_{i-2}[(h_i - 1) + (h_{i-1} - 1)\hat{t}'_{i-3} \prod_{l=1}^{i-2} p_l + \sum_{j=2}^{i-2} [(h_j - 1)\hat{t}'_{j-2} \prod_{l=1}^{j-1} p_l \prod_{q=j+1}^{i-1} p_q \hat{x}'_q]] + t_{i-1} p_i \quad (16)$$

.
.

$$x_r = \hat{x}'_r [(h_r - 1) + (h_{r-1} - 1)\hat{t}'_{r-3} \prod_{l=1}^{r-2} p_l + \sum_{j=2}^{r-2} [(h_j - 1)\hat{t}'_{j-2} \prod_{l=1}^{j-1} p_l \prod_{q=j+1}^{r-1} p_q \hat{x}'_q]] + t_{r-1} \prod_{j=1}^{r-1} p_j \quad (17)$$

$$t_{r-2} = \hat{t}'_{r-2}[(h_r - 1) + (h_{r-1} - 1)\hat{t}'_{r-3} \prod_{l=1}^{r-2} p_l + \sum_{j=2}^{r-2} [(h_j - 1)\hat{t}'_{j-2} \prod_{l=1}^{j-1} p_l \prod_{q=j+1}^{r-1} p_q \hat{x}'_q]] + t_{r-1} p_r \quad (18)$$

where $\hat{x}'_i$ and $\hat{t}'_{i-2}$ are some given solution of

$$p_i x'_i - t'_{i-2} \prod_{k=1}^{i-1} p_k = 1 \qquad (19)$$

and $j = 2,3,4,...,i-2$, $i = 4,5,6,...,r$ and supposing $\hat{t}'_0 = \hat{x}'_1$. It is clear that the given solutions $\hat{x}'_i$ and $\hat{t}'_{i-2}$ always exist, as equation (19) is a special case of equation (9). Note that in (17) and (18), $t_{r-1}$ is a final free integer parameter. Using (17) and (18) the value $t_i$ can be re-written in term of $t_{r-1}$. Moreover, using (16) the variable $x_i$ can be re-written in terms of $t_{r-1}$ and the general solutions of (19) and $h_i$ and $p_i$. Thus $x_1$ can be obtained as:

$$x_1 = t_{r-1} \prod_{l=2}^{r} p_l + (h_r - 1)\hat{t}'_{r-2} \prod_{l=2}^{r-1} p_l + (h_{r-1} - 1)\hat{t}'_{r-3} p_r \hat{x}'_r \prod_{l=2}^{r-2} p_l + \sum_{j=2}^{r-2} [(h_j - 1)\hat{t}'_{j-2} \prod_{l=1}^{j-1} p_l \prod_{q=j+1}^{r} p_q \hat{x}'_q] \qquad (20)$$

Using (20), (19) and (7), linear function $Z_r(t_{r-1}, h_j)$ (as the general functional form of set $H_r$) can be obtained as:

$$Z_r(t_{r-1}, h_i) = t_{r-1} \prod_{l=1}^{r} p_l + (p_r \hat{x}'_r - 1)(h_r - 1) + \sum_{j=2}^{r-1} [(p_j \hat{x}'_j - 1)(h_j - 1) \prod_{q=j+1}^{r} p_q \hat{x}'_q] - 1 \qquad (21)$$

Equation (19) and formula (21) are equivalent to formulas (A-1) and (A-2), thus proof of Theorem 1 is completed ∎.

## 3. Some Corollaries and Remarks

**Remark 1.** $Z_r(t, h_i)$ **algebraically defines all integer numbers which are not divisible by primes:** $p_1, p_2, ..., p_r$.

We just recall that all terms in formula (A-2) are only made up of presupposed prime numbers $p_1, p_2, p_3, ..., p_r$ (values of $x'_i$ are depended on them), and integer variables $t$ and $h_i$. Furthermore, $x'_i$ in formula (A-2) doesn't have unique value (because in principle, equation (A-1) has infinite integer solutions), hence function $Z_r(t, h_i)$ can be written in different but equivalent forms.

**Corollary 1.** Suppose $q_1, q_2, ..., q_r$ are some given co-prime numbers, and let $x'_i$ and $t'_i$ are some integer solution of linear equation:

$$q_i x_i - t_i \prod_{k=1}^{i-1} q_k = 1 \qquad (22)$$

where $i = 2, 3, 4, ..., r$, **then the linear function** $W_r(t, h_i)$:

$$W_r(t, h_i) = t \prod_{l=1}^{r} q_l + (q_r x'_r - 1)(h_r - h_1) + \sum_{j=2}^{r-1}[(q_j x'_j - 1)(h_j - h_1)\prod_{f=j+1}^{r} q_f x'_f] - h_1 \qquad (23)$$

**defines all (and only) integer numbers which are not divisible by co-primes:** $q_1, q_2, ..., q_r$; where $t$ is an integer variable: $t \in \mathbb{Z}$ and $h_i \in \{1, 2, 3, 4, ..., q_i - 1\}$ and $r = 3, 4, 5, ...$; for $r = 1, 2$ we have

$$W_1(t, h_i) = q_1 t - h_1, \quad W_2(t, h_i) = q_1 q_2 t + (q_2 x'_2 - 1)(h_2 - h_1) - h_1 \qquad (24)$$

**Corollary 2.** Using some theorems such as Bertrand's postulate (that states for every n > 1, there is always at least one prime $p$ such that n < p < 2n), the action range of (A-3) can be expanded for larger intervals. Here by Bertrand's postulate, it is easy to show that all (and only) k-almost prime numbers (where $\exists n \in \mathbb{N}$, $1 \leq k \leq n$) can be redefined by $Z_r(t, h_i)$, in interval $[p_{r+s}^n, p_{r+s}^{n+1})$, if $p_{r+1} > 2^{(n+1)(s-1)}$; where $s \geq 1$, $i_j = 1, 2, 3, ..., s$ and $j = 1, 2, 3, ..., n$.

**Corollary 3.** Comparing function $Z_r(t, h_i)$ with function $v_i = f^{-1}(u_i)$ in [40], we get the following new relations for prime numbers $p_1, p_2, p_3, ..., p_r$ ($e = 2, 3, 4, ..., r-1$),

$\exists x'_2, x'_3, ..., x'_r$ (where $x'_i$ is some integer solution of equation (A-1): $p_i x_i - t_i \prod_{l=1}^{i-1} p_l = 1$), $\exists s \in \mathbb{Z}$:

$$(2s-1)\prod_{i=2}^{r-1} p_i = x'_r [(p_{r-1} x'_{r-1} - 1) + \sum_{i=2}^{r-2}[(p_i x'_i - 1)\prod_{j=i+1}^{r-1} p_j x'_j]] \qquad (25)$$

$$(\prod_{l=1}^{r} p_l / p_e)[(\prod_{l=1}^{r} p_l / p_e) \text{ modInverse } (p_e)] = -(p_e x'_e - 1)\prod_{q=e+1}^{r} p_q x'_q \qquad (26)$$

These relations could be simplified, too.

**Corollary 4.** The number of primes in interval $(1, p_{r+1}^2)$, using function $Z_r(t, h_i)$. As function $Z_r(t, h_i)$ is linear, from formula (A-1), the number of the members of set $H_r$ (formula (1)) in intervals $(1, \prod_{l=1}^{r} p_l)$, $(\prod_{l=1}^{r} p_l + 1, 2\prod_{l=1}^{r} p_r)$, $(2\prod_{l=1}^{r} p_l + 1, 3\prod_{l=1}^{r} p_r)$ and so on, is:

$$\prod_{l=1}^{r}(p_l - 1) - 1 \qquad (27)$$

Now using (27), approximately, the number of primes in interval $(1, p_{r+1}^2)$ is:

$$\pi(p_{r+1}^2) \cong r + p_{r+1}^2 [\prod_{l=1}^{r}(p_l - 1) - 1] / \prod_{l=1}^{r} p_l \qquad (28)$$

**Remark 2.** The values of $Z_r(t, h_i)$. Using formula (A-2) or (25) and (26), we may easily write (in an unique way here) the values of function $Z_r(t, h_i)$ for $r = $ 1, 2, 3, 4, 5, ..., in certain forms:

$Z_1(t, h_i) = 2t + 1$
$Z_2(t, h_i) = 6t + 4h_2 + 3$
$Z_3(t, h_i) = 30t + 6h_3 + 10h_2 + 15$
$Z_4(t, h_i) = 210t + 120h_4 + 126h_3 + 70h_2 + 105$
$Z_5(t, h_i) = 2310t + 210h_5 + 330h_4 + 1386h_3 + 1540h_2 + 1155$
$Z_6(t, h_i) = 30030t + 6930h_6 + 16380h_5 + 25740h_4 + 6006h_3 + 20020h_2 + 15015$
$Z_7(t, h_i) = 510510t + 450450h_7 + 157080h_6 + 46410h_5 + 145860h_4 + 306306h_3 + 170170h_2 + 255255$
$Z_8(t, h_i) = 9699690t + 9189180h_8 + 9129120h_7 + 3730650h_6 + 6172530h_5 + 8314020h_4 +$
$\qquad + 3879876h_3 + 323230h_2 + 4849845$
…

$\qquad (29)$

Recently, some standard integer sequences (sequence A240673 in The On-Line Encyclopedia of Integer Sequences® [42]) are constructed based on the coefficients of function $Z_r(t, h_i)$.

# Acknowledgment

Special thanks are extended to Prof. and Academician Andrzej Schinzel (Institute of Mathematics of the Polish Academy of Sciences), Prof. Carl Pomerance (Dartmouth College), Prof. Graham Jameson (Lancaster University), Prof. Maurice H.P.M. van Putten (MIT), Prof. Hans Juergen Weber (University Of Virginia), Prof. Len Adleman (University of Southern California), Prof. Gérald Tenenbaum (University of Lorraine) for their interest to this work, as well as their valuable guidance and feedbacks.